\documentclass{compositio}
\usepackage{graphicx}
\usepackage{amsfonts}
\usepackage{eucal}
\usepackage{amscd}
\usepackage{amssymb}
\usepackage[all]{xy}
\newtheorem{thm}{Theorem}[section]
\newtheorem{cor}[thm]{Corollary}
\newtheorem{lem}[thm]{Lemma}
\newtheorem{prop}[thm]{Proposition}
\theoremstyle{definition}
\newtheorem{defn}[thm]{Definition}
\theoremstyle{remark}

\newcommand{\norm}[1]{\left\Vert#1\right\Vert}

\newcommand{\Lie}[1]{ {\mathfrak #1} }

\newcommand{\variety}[1]{\mathbf{#1}}

\newcommand{\nats}{\mathbb N}

\newcommand{\rats}{\mathbb Q}
\newcommand{\reals}{\mathbb R}
\newcommand{\complex}{\mathbb C}

\newcommand{\A}{\mathcal{A}}

\newcommand{\M}{\mathcal{M}}

\newcommand{\U}{\mathcal{U}}
\newcommand{\OO}{\mathcal{O}}
\newcommand{\ZU}{\mathcal{Z}}

\newcommand{\isom}{\cong}
\newcommand{\gK}{{( \Lie{g}, K)}}
\newcommand{\J}{\mathcal J}
\newcommand{\ii}{\mathbf{i}}
\newcommand{\pp}{\mathbf{p}}
\newcommand{\ev}{\mathbf{ev}}
\newcommand{\Ell}{\mathbf{l}}
\newcommand{\dd}{\mathbf{d}}
\newcommand{\II}{\mathbf{I}}
\newcommand{\PP}{\mathbf{P}}

\newcommand{\LL}{\mathbf{L}}
\newcommand{\DD}{\mathbf{D}}
\begin{document}

\title{Multiplying Modular Forms}%
\author{Martin H. Weissman}%
\email{weissman@ucsc.edu}%
\address{Department of Mathematics, University of California, Santa Cruz, CA 95064}%


\begin{abstract}
The space of elliptic modular forms of fixed weight and level can be identified with a space of intertwining operators, from a holomorphic discrete series representation of $SL_2(\reals)$ to a space of automorphic forms.  Moreover, multiplying elliptic modular forms corresponds to a branching problem involving tensor products of holomorphic discrete series representations.

In this paper, we explicitly connect the ring structure on spaces of modular forms with branching problems in the representation theory of real semisimple Lie groups.  Furthermore, we construct a family of intertwining operators from discrete series representations into tensor products of other discrete series representations.  This collection of intertwining operators provides the well-known ring structure on spaces of holomorphic modular forms.

An analytic subtlety prevents this collection of intertwining operators from directly yielding a ring structure on some spaces of modular forms.  We discuss discrete decomposability, and its role in constructing ring structures.
\end{abstract}

\maketitle
\section*{Introduction}

One of the easiest ways to construct a cuspidal elliptic modular form is by multiplying and adding Eisenstein series.  For example, a cusp form of weight $12$ is given by $\Delta = g_2^3 - 27 g_3^2$, where $g_2$ and $g_3$ are Eisenstein series of weights $4$ and $6$, respectively.  The fact that elliptic modular forms can be multiplied is practically obvious -- but it relies on the fact that the product of holomorphic functions is again holomorphic.  Multiplication of modular forms extends, without difficulty, to multiplying Siegel modular forms, and more generally, to multiplying scalar-valued holomorphic modular forms on tube domains.

In Remarque 2.1.4 of \cite{Del}, Deligne notices a connection between multiplication of elliptic modular forms and the tensor product of two holomorphic discrete series representations.  Motivated by this observation, we approach the problem of multiplying modular forms as a branching problem for tensor products of discrete series representations.  In expanding on Deligne's remark, we attempt to answer two questions:
\begin{enumerate}
\item
What is the precise connection between branching problems in representation theory and multiplication of modular forms?  Specifically, what statements from representation theory allow the construction of a ring of modular forms?
\item
How, and when, can one find a canonical embedding of one discrete series representation into the tensor product of two other discrete series representations of a real semisimple Lie group?
\end{enumerate}

It is far from clear whether one can multiply non-holomorphic modular forms.  For example, one might try to multiply those modular forms on the split group of type $G_2$ considered by Gross, Gan, and Savin in \cite{GGS}, or those on the split group of type $D_4$ in \cite{Wei}.  Most generally, we are interested in multiplying modular forms, whose associated representation of a real semisimple group is in the discrete series.  The ``algebraic modular forms'' of Gross \cite{Gr1}, where the associated real semisimple group is compact, can be multiplied by results of Khuri-Makdisi \cite{KMK}.  This paper can be seen as a generalization of his results, from finite-dimensional representations of compact Lie groups, to holomorphic discrete series representations of semisimple Lie groups, and an attempt to generalize to other discrete series.

One of our main results is a canonical embedding of a discrete series representation in the tensor product of two others.  The existence of an embedding follows directly from work of Vargas \cite{Var}; using the analytic results of \O rsted-Vargas \cite{O-V}, we are able to make such an embedding canonical.

We finish this paper by describing how our representation-theoretic multiplication corresponds to more classical multiplication of modular forms.  This expands on Deligne's remark.

\begin{acknowledgements}
I would like to thank many people for their advice during the preparation of this paper.  I exchanged a series of e-mails with J. Vargas while working on these results -- I am very grateful for his advice.  It also appears that the current work of Vargas has some overlap with mine, though I believe that this paper is of significant independent interest.  I would also like to thank W. Schmid, H.-Y. Loke, W.-T. Gan, B. Gross, and D. Toledo, G. van der Geer, and L. Ji for some additional advice along the way.

The original draft of this paper contained a serious flaw; I am very thankful to Gordan Savin for pointing out this problem.  Although this flaw curtails the inital ambitions of this paper (multiplying certain quaternionic modular forms), it is useful to present some negative results.  In particular, results of Kobayashi \cite{Kob} on discrete decomposability explain why some modular forms cannot be multiplied.
\end{acknowledgements}
\section*{Notation}

$\variety{G}$ will always denote a connected semisimple algebraic group over $\rats$.  Then we let $\Lie{g}$ denote its complexified Lie algebra:  $\Lie{g} = Lie(\variety{G}) \otimes_\rats \complex$.  Let $\U$ denote the universal enveloping algebra of $\Lie{g}$, and $\ZU$ the center of $\U$.

We write $G$ for $\variety{G}(\reals)$, viewed as a semisimple real Lie group.  Note that $G$ is connected with finite center.  Let $K$ denote a maximal compact subgroup of $G$; $K$ is then a connected compact Lie group.  $\Gamma$ will denote an arithmetic subgroup of $\variety{G}(\rats)$.

All vector spaces will be complex vector spaces; all representations will be on complex vector spaces, and all tensor products will be over $\complex$.  If $H$ is a group, and $V,V'$ are representations of $H$, then we write $V \boxtimes V'$ for the ``external tensor product'' representation of $H \times H$.  On the other hand, we write $V \otimes V'$ for the usual tensor product of representations, i.e.,  the restriction of $V \boxtimes V'$ to the diagonally embedded $H$ in $H \times H$.

If $U,V,W$ are (complex) vector spaces, then there are canonical isomorphisms which we name:
\begin{eqnarray*}
can & \colon & U \rightarrow \complex \otimes U, \\
comm & \colon & U \otimes V \rightarrow V \otimes U, \\
assoc & \colon & (U \otimes V) \otimes W \rightarrow U \otimes (V \otimes W).
\end{eqnarray*}

\section{Automorphic and Modular Forms}

The definition of an automorphic form varies slightly from one source to another.  We follow the definition contained in the article of Borel and Jacquet \cite{B-J}.  Our definition of modular forms is based on that of Gross, Gan, and Savin \cite{GGS}.

\subsection{$\gK$-modules}

The definition of a $\gK$-module is taken from Borel-Wallach \cite{B-W}.
\begin{defn}
A $\gK$-module is a complex vector space $W$, endowed with simultaneous actions of $\Lie{g}$ and $K$ (meaning a Lie algebra homomorphism $d \pi \colon \Lie{g} \rightarrow End_\complex(W)$ and a group homomorphism $\pi_K \colon K \rightarrow GL(W)$), satisfying the following axioms
\begin{itemize}
\item
$W$ is locally finite and semisimple as a $K$-module.
\item
For all $X \in \Lie{g}$, all $k \in K$, and all $w \in W$,
$$\pi_K(k) (d \pi(X) w) = d \pi [Ad(k)(X)] (\pi_K(k) w).$$
\item
If $W' \subset W$ is a $K$-stable finite-dimensional subspace of $W$, then the representation of $K$ on $W'$ is differentiable, and the resulting action $d (\pi_K)$ of $\Lie{k}$ agrees with the restriction $(d \pi) \vert_\Lie{k}$ of $d \pi$.
\end{itemize}
\end{defn}

If $W$ is a $\gK$-module, then the Lie algebra action of $\Lie{g}$ extends canonically to an action of the universal enveloping algebra $\U$.  We say that $W$ has an {\em infinitesimal character} if there exists a homomorphism $\omega \colon \ZU \rightarrow \complex$ such that $Z w = \omega(Z) \cdot w$, for all $Z \in \ZU$, and all $w \in W$.

\begin{defn}
A Harish-Chandra module is a $\gK$-module $W$, such that $W$ is finitely generated as a $\U$-module, and such that every irreducible representation of $K$ occurs with finite multiplicity in $W$.
\end{defn}

If $W$ is an irreducible Harish-Chandra module, then $W$ has an infinitesimal character.
\begin{defn}
A unitary structure on a $\gK$-module $W$ is a positive non-degenerate scalar product $\langle \cdot, \cdot \rangle$ on $W$, which satisfies:
\begin{itemize}
\item
For all $k \in K$, and all $v,w \in W$, $\langle \pi_K(k) v, \pi_K(k) w \rangle = \langle v, w \rangle$.
\item
For all $X \in \Lie{g}$, and all $v,w \in W$, $\langle d \pi(X) v, w \rangle + \langle v, d \pi(X) w \rangle = 0$.
\end{itemize}
\end{defn}

\begin{defn}
A {\em weight} is an irreducible Harish-Chandra module $W$, endowed with a unitary structure.
\end{defn}

As tensor products of $\gK$-modules arise often, we fix a few definitions here.  First, there exists a unique associative algebra homomorphism $\Delta \colon \U \rightarrow \U \otimes \U$, which satisfies $\Delta(X) = 1 \otimes X + X \otimes 1$, for all $X \in \Lie{g}$.  Thus, if $W_1$ and $W_2$ are $\gK$-modules, the tensor product $W_1 \otimes W_2$ is endowed with a natural $\U$ action as well as a $K$ action via this diagonal embedding.  This makes the tensor product $W_1 \otimes W_2$ a $\gK$ module as well.

\subsection{Automorphic and Modular Forms}

Here, we recall the definition of automorphic forms given by Borel-Jacquet in \cite{B-J}.  Fix a finite-dimensional complex representation $\sigma \colon G \rightarrow GL(E)$, with finite kernel and closed image.  Fix a $\sigma(K)$-invariant Hilbert space structure on $E$.  Define, for all $g \in G$,
$$\norm{g} = \sqrt{ \vert tr(\sigma(g)^\ast \sigma(g) ) \vert }.$$
Here $\sigma(g)^\ast$ denotes the adjoint to $\sigma(g)$, with respect to the chosen Hilbert space structure.

\begin{defn}
An {\em automorphic form} on $G$ (with respect to $\Gamma$) is a complex-valued function $f \in C^\infty(\Gamma \backslash G)$ such that:
\begin{enumerate}
\item
$f$ has moderate growth.  In other words, there exists a positive integer $n$, and a constant $C$ such that:
$$\vert f(g) \vert \leq C \cdot \norm{g}^n, \mbox{ for all } g \in G.$$
\item
$f$ is $K$-finite on the right.
\item
$f$ is $\ZU$-finite on the right.
\end{enumerate}
\end{defn}

We write $\A(\Gamma \backslash G)$ for the complex vector space of automorphic forms on $\Gamma \backslash G$.  This space is a $\gK$-module, by right-translation.  Since it will be convenient later, we also define:
\begin{defn}
A {\em relaxed automorphic form} on $G$ (wish respect to $\Gamma$) is a complex-valued function $f \in C^\infty(\Gamma \backslash G)$ such that conditions (1) and (2), described above, are satisfied.
\end{defn}
We write $\A_{rel}(\Gamma \backslash G)$ for the complex vector space of relaxed automorphic forms on $\Gamma \backslash G$.  This space is also a $\gK$-module, by right-translation.

Modular forms arise from irreducible $\gK$-submodules of $\A(\Gamma \backslash G)$.  Let $W$ be a weight, i.e., an irreducible unitary Harish-Chandra module.
\begin{defn}
A modular form on $\Gamma \backslash G$, of weight $W$, is a homomorphism $F$ of $\gK$-modules from $W$ to $\A(\Gamma \backslash G)$.
\end{defn}
The space of modular forms of weight $W$ is given by:
$$\M(\Gamma \backslash G, W) = Hom_\gK(W, \A(\Gamma \backslash G)).$$

We may also define relaxed modular forms:
\begin{defn}
A relaxed modular form on $\Gamma \backslash G$, of weight $W$, is a homomorphism $F$ of $\gK$-modules from $W$ to $\A_{rel}(\Gamma \backslash G)$.
\end{defn}
The space of relaxed modular forms of weight $W$ is given by:
$$\M_{rel}(\Gamma \backslash G, W) = Hom_{\gK}(W, \A_{rel}(\Gamma \backslash G)).$$
However, it turns out that relaxed modular forms are no different than ordinary modular forms:
\begin{prop}
Every relaxed modular form is a modular form.  In other words, $\M_{rel}(\Gamma \backslash G, W) = \M(\Gamma \backslash G, W)$.
\end{prop}
\begin{proof}
Suppose that $F \in \M_{rel}(\Gamma \backslash G, W)$.  We must show that, for all $w \in W$, the relaxed automorphic form $F(w)$ is $\ZU$-finite.  Since $W$ is a weight, it has an infinitesimal character $\omega$.  Since $F$ is a $\gK$-module homomorphism, it is $\ZU$-intertwining.  Therefore, for all $Z \in \ZU$, we have:
$$Z F(w) = F(Z w) = \omega(Z) \cdot F(w).$$
Hence $F(w)$ is trivially $\ZU$-finite.
\end{proof}

\section{Multiplication}
Hereafter, we omit mention of $\Gamma$ and $G$, when it does not reduce clarity.

Suppose that $W_1$ and $W_2$ are weights (irreducible unitary Harish-Chandra modules).  We are interested in the following question:  under what circumstances does there exist a weight $W$, and a natural $\complex$-bilinear operation:
$$\M(W_1) \otimes \M(W_2) \rightarrow \M(W)?$$

\subsection{Multiplying automorphic forms}

If $f_1, f_2 \in \A = \A(\Gamma \backslash G)$, let $m(f_1, f_2) = f_1 \cdot f_2$ denote the pointwise product.  This extends to a bilinear map:
$$m \colon \A \otimes \A \rightarrow C^\infty(\Gamma \backslash G).$$
\begin{prop}
The map $m$ is a homomorphism of $\gK$-modules.  The image of $m$ is contained in the space $\A_{rel}$ of relaxed automorphic forms.
\end{prop}
\begin{proof}
The bilinear form $m$ is a morphism of $\gK$-modules, since multiplication is compatible with right-translation.  In other words, translation on the right by $K$ is preserved by multiplication.  Furthermore, if $X \in \Lie{g}$, then we have:
$$X (f_1 \cdot f_2) = (X f_1) \cdot f_2 + f_1 \cdot (X f_2) = m \left( (X \otimes 1 + 1 \otimes X) (f_1 \otimes f_2) \right).$$

Now, we check that if $f_1, f_2 \in \A$, then $f = f_1 \cdot f_2 \in \A_{rel}$.  Since $f_1, f_2$ have moderate growth, we have:
$$\vert f(g) \vert = \vert f_1(g) \vert \vert f_2(g) \vert \leq C_1 C_2 \cdot \norm{g}^{n_1 + n_2},$$
where $C_1, C_2, n_1, n_2$ are the appropriate constants of moderate growth for $f_1, f_2$.  It follows that $f$ has moderate growth.

Finally, let $\Phi_1, \Phi_2, \Phi$ denote the complex vector spaces spanned by the $K$-translates of $f_1$, $f_2$, and $f$, respectively.  Then $\Phi \subset m(\Phi_1 \otimes \Phi_2)$, and $\Phi_1$ and $\Phi_2$ are finite-dimensional.  Hence $\Phi$ is finite-dimensional.
\end{proof}

The above proposition states that the product of two automorphic forms is {\em almost} an automorphic form; the only condition missing is $\ZU$-finiteness.

\subsection{Multiplying modular forms}

Suppose that we are given three weights $W, W_1, W_2$, and a homomorphism of $\gK$-modules:
$$\mu_{1,2} \colon W \rightarrow W_1 \otimes W_2.$$
Then we may define a bilinear map:
$$m_{1,2} = Co(\mu_{1,2}) \colon \M(W_1) \otimes \M(W_2) \rightarrow \M_{rel}(W),$$
by $m_{1,2}(F_1, F_2) = F$, where
$$F = m \circ [F_1 \otimes F_2] \circ \mu_{1,2},$$
following the sequence of maps:
$$\xymatrix@C=0.5cm{
  W \ar[rr]^{\mu_{1,2}} \ar@/_1pc/[rrrrrr]|{F} && W_1 \otimes W_2 \ar[rr]^{F_1 \otimes F_2} && \A \otimes \A \ar[rr]^{m} && \A_{rel} }
$$
Indeed, $F$ is a $\gK$-module homomorphism, as a composite of three $\gK$-module homomorphisms.  Since $\M_{rel}(W) = \M(W)$, we arrive at:
\begin{prop}
If $F_1 \in \M(W_1)$, and $F_2 \in \M(W_2)$, then $F = m_{1,2}(F_1, F_2) \in \M(W)$.  In this way, we get a linear map:
$$Co \colon Hom_{\gK}(W, W_1 \otimes W_2) \rightarrow Hom_\complex(\M(W_1) \otimes \M(W_2), \M(W)).$$.
\end{prop}

\subsection{A graded ring}
Suppose that $C$ is a commutative monoid, with unit element $0$.  If $\lambda_1, \lambda_2 \in C$, we write $\lambda_1 + \lambda_2$ for the composition.  A $C$-graded ring (over $\complex$), is a commutative unital $\complex$-algebra $R$, together with an isomorphism of vector spaces:
$$R = \bigoplus_{\lambda \in C} R_\lambda,$$
for subspaces $R_\lambda$ of $R$, satisfying:
$$\forall \lambda_1, \lambda_2 \in C, R_{\lambda_1} \cdot R_{\lambda_2} \subset R_{\lambda_1 + \lambda_2}.$$
We also assume that $R_0 = \complex$.  To construct such an algebra, it suffices to define the collection of vector spaces $R_\lambda$, and the bilinear multiplication maps among them.

In order to construct a $C$-graded ring of modular forms, we require the following:  suppose that for every $\lambda \in C$, one has an weight $W_\lambda$.  Suppose that $W_0 = \complex$ is the trivial Harish-Chandra module.  Moreover, suppose we are given, for every pair $\lambda_1, \lambda_2$, a homomorphism of $\gK$-modules:
$$\mu_{\lambda_1, \lambda_2} \colon W_{\lambda_1 + \lambda_2} \rightarrow W_{\lambda_1} \otimes W_{\lambda_2}.$$
We refer to the set $\{ \mu_{\lambda_1, \lambda_2} \}_{\lambda_1, \lambda_2 \in C}$ as the {\em family of comultiplications}.  Define the graded vector space $\M_C = \bigoplus_{\lambda \in C} \M(W_\lambda)$.  Then we have a set of multiplications on $\M_C$ via the bilinear maps $m_{\lambda_1, \lambda_2} = Co(\mu_{\lambda_1, \lambda_2})$ on each graded piece of $\M$.  This makes $\M_C$ a $C$-graded ring if and only if the family of comultiplications satisfy the following axioms:
\begin{description}
\item[Co-identity]
For all $\lambda$, there is an equality of maps:
$$
\xymatrix{
W_\lambda  \ar@/^/[rr]|{\mu_{0,\lambda}} \ar@/_/[rr]|{can} && \complex \otimes W_\lambda
}
$$
\item[Co-commutativity]
For all $\lambda_1, \lambda_2 \in C$, and $\lambda = \lambda_1 + \lambda_2$, the following diagram commutes:
$$
\xymatrix@R=0.5cm{
                &     W_{\lambda_1} \otimes W_{\lambda_2}     \ar@{<->}[dd]^{comm}     \\
  W_\lambda \ar[ur]^{\mu_{\lambda_1, \lambda_2}} \ar[dr]_{\mu_{\lambda_2, \lambda_1}}                 \\
                &     W_{\lambda_2} \otimes W_{\lambda_1}                   }
$$
\item[Co-associativity]
For all $\lambda_1, \lambda_2, \lambda_3 \in C$, the following diagram commutes:
$$
\xymatrix@C+1cm{
W_{\lambda_1 + \lambda_2 + \lambda_3} \ar @{=} [d] \ar[r]^{\mu_{\lambda_1, \lambda_2+\lambda_3}}
& W_{\lambda_1} \otimes W_{\lambda_2+\lambda_3} \ar[r]^{Id \otimes \mu_{\lambda_2, \lambda_3}}
& W_{\lambda_1} \otimes (W_{\lambda_2} \otimes W_{\lambda_3}) \ar@{<->}[d]^{assoc}
\\
W_{\lambda_1 + \lambda_2 + \lambda_3} \ar[r]^{\mu_{\lambda_1+\lambda_2, \lambda_3}}
& W_{\lambda_1+\lambda_2} \otimes W_{\lambda_3} \ar[r]^{\mu_{\lambda_1, \lambda_2} \otimes Id}
& (W_{\lambda_1} \otimes W_{\lambda_2}) \otimes W_{\lambda_3}
}
$$
\end{description}

\section{Compact Groups}

Our first task will be to consider tensor products of representations of $K$, a maximal compact subgroup of $G$.  In this section, we use the Borel-Weil realization of representations of $K$, in order to construct canonical embeddings of some representations of $K$ in tensor products of others.  While essentially contained in \cite{KMK}, we review these ideas thoroughly here as it will be important in what follows.

\subsection{Representations of $K$}
Recall that $G$ is a connected semisimple Lie group, with finite center, and $K$ is a maximal compact subgroup of $G$.  In particular, $K$ is also connected.  Let $T$ be a maximal torus in $K$.  Assume hereafter that $G$ has the same rank as $K$, so $T$ is also a maximal torus in $G$.

We write $\Lambda$ for the abelian group $Hom_{alg}(T_\complex, \complex^\times)$.  Fix a Borel subgroup $Q_\complex$ of $K_\complex$ containing $T_\complex$, with unipotent racical $U_\complex$.  Let $U_\complex^-$ denote the unipotent radical of the opposite Borel subgroup.

Let $\Delta_k^+$ denote the associated set of positive roots, for the adjoint action of $T_\complex$ on $\Lie{q}_\complex$.  Let $\Lambda_k^+$ denote the associated cone of dominant weights in $\Lambda$:
$$\Lambda_k^+ = \{ \lambda \in \Lambda \mbox{ such that } \forall \delta \in \Delta_k^+, \langle \lambda, \delta \rangle \geq 0 \}.$$

We recall the famous Borel-Weil theorem that realizes the irreducible representations of $K$ in the cohomology of a flag variety.  This realization is particularly important, since we need some compatibilities among different representations.

The compact quotient $X = K/T$ may be identified with $K_\complex / Q_\complex$, and may thus be considered as a complex projective variety.  Suppose that $\lambda \in \Lambda$.  Extend $\lambda$ to a homomorphism from $Q_\complex$ to $\complex^\times$, by having $U_\complex$ sent to $1$.

Let $\OO(-\lambda)$ denote the line bundle on $X$ associated to the character $-\lambda$.  The holomorphic sections of $\OO(-\lambda)$ consist of those holomorphic functions $h \colon K_\complex \rightarrow \complex$ satisfying:
$$h(kq) = \lambda(q) \cdot h(k) \mbox{ for all } k \in K_\complex, q \in Q_\complex.$$
We write $H^0(X, \OO(-\lambda))$ for this space of holomorphic sections.

The complex vector space $H^0(X, \OO(-\lambda))$ is a finite-dimensional representation of $K$, where $K$ acts by left-translation $\tau = \tau_\lambda$:
$$[\tau(k) f](x) = f(k^{-1} x).$$

We recall the following theorem, which we believe originally appeared in \cite{Ser}; we refer to Knapp's book \cite{Kna} for a treatment as well:
\begin{thm}[(Borel-Weil)]
Suppose that $\lambda \in \Lambda_k^+$.  Then $H_\lambda = H^0(X, \OO(-\lambda))$ is an irreducible representation of $K$ with highest weight $\lambda$.
\end{thm}

The Borel-Weil realization of the irreducible representations of $K$ also fixes a number of auxiliary structures.  First, the representations $H_\lambda$ are endowed with natural Hilbert space structures.  Namely, if $h_1, h_2 \in H_\lambda$, then we may define:
$$\langle h_1, h_2 \rangle_\lambda = \int_K h_1(k) \overline{h_2(k)} dk.$$
Here we normalize Haar measure so that $K$ has measure $1$.  Note that this inner product is $K$-invariant.

The $\lambda$ weight space of $H_\lambda$ is one-dimensional.  Moreover, the Borel-Weil realization yields two natural choices of highest-weight vectors.

\begin{prop}
There exists a unique highest-weight vector $h_\lambda \in H_\lambda$, which satisfies $h_\lambda(1) = 1$.
\end{prop}
\begin{proof}
Any highest-weight vector $h_\lambda$ not only satisfies $h_\lambda(kq) = \lambda(q) h(k)$ for all $q \in Q_\complex$, but also $h_\lambda(u^{-}k) = h_\lambda(k)$, for all $u^-$ in $U_\complex^-$.  Thus the values of a highest-weight vector on the open Bruhat cell $U_\complex^- Q_\complex$ are nonzero complex multiples of the value at $1$.  It follows that any nonzero highest-weight vector $h_\lambda$ satisfies $h_\lambda(1) \neq 0$.  Since the highest-weight space is one-dimensional, there is a unique highest-weight vector $h_\lambda \in H_\lambda$, which additionally satisfies $h_\lambda(1) = 1$.
\end{proof}

\begin{prop}
There exists a unique highest-weight vector $h_\lambda^\vee \in H_\lambda$, which satisfies $\langle h_\lambda^\vee, h \rangle = h(1)$, for all $h \in H_\lambda$.
\end{prop}
\begin{proof}
The map $h \mapsto h(1)$ is a linear functional on the finite-dimensional Hilbert space $H_\lambda$.  Hence there exists a unique $h_\lambda^\vee \in H_\lambda$ which satisfies:
$$\langle h_\lambda^\vee, h \rangle = h(1), \mbox{ for all } h \in H_\lambda.$$
Moreover, the $K$-invariance of the inner product yields:
$$\langle \tau(t) h_\lambda^\vee, h \rangle = \langle h_\lambda^\vee, \tau(t)^{-1} h \rangle = h(t) = \lambda(t) h(1),$$
for all $t \in T$.  Hence $h_\lambda^\vee$ is a highest-weight vector in $H_\lambda$.
\end{proof}

Note that $\langle h_\lambda^\vee, h_\lambda \rangle = 1$.  We can also compute the following:
\begin{prop}
Let $d_\lambda = dim(H_\lambda)$.  Then we have:
$$\langle h_\lambda, h_\lambda \rangle = d_\lambda^{-1}, \mbox{ and } \langle h_\lambda^\vee, h_\lambda^\vee \rangle = d_\lambda.$$
\label{HWF}
\end{prop}
\begin{proof}
We may use the Schur orthogonality relations to compute this inner product:
\begin{eqnarray*}
\langle h_\lambda, h_\lambda \rangle & = & \int_K h_\lambda(k) \overline{h_\lambda(k)} dk, \\
& = & \int_K \langle h_\lambda^\vee, \tau(k)^{-1} h_\lambda \rangle \overline{ \langle h_\lambda^\vee, \tau(k)^{-1} h_\lambda \rangle } dk, \\
& = & \int_K \langle  \tau(k) h_\lambda^\vee, h_\lambda \rangle \overline{ \langle \tau(k) h_\lambda^\vee,  h_\lambda \rangle } dk, \\
& = & { 1 \over {d_\lambda} } \langle h_\lambda^\vee, h_\lambda^\vee \rangle \overline{ \langle h_\lambda, h_\lambda \rangle }.
\end{eqnarray*}
Since $\langle h_\lambda, h_\lambda \rangle$ is real and nonzero, we arrive at:
$$\langle h_\lambda^\vee, h_\lambda^\vee \rangle = d_\lambda.$$
Since $h_\lambda$ is a scalar multiple of $h_\lambda^\vee$, and $\langle h_\lambda, h_\lambda^\vee \rangle = 1$, we deduce:
$$\langle h_\lambda, h_\lambda \rangle = d_\lambda^{-1}.$$
\end{proof}

\subsection{Tensor products of Borel-Weil representations}

Consider the tensor product $H_{\lambda_1} \otimes H_{\lambda_2}$, for two dominant weights $\lambda_1, \lambda_2$.  Let $\lambda = \lambda_1 + \lambda_2$.  Then we find a natural map:
$$\mu_{\lambda_1, \lambda_2}^\ast \colon H_{\lambda_1} \otimes H_{\lambda_2} \rightarrow H_\lambda.$$
This map is defined by extending the following linearly:
$$\mu_{\lambda_1, \lambda_2}^\ast(h_1 \otimes h_2) = h_1 \cdot h_2.$$
This map is $K$-intertwining and nonzero.  Moreover, evaluating at the identity yields:
$$\mu_{\lambda_1, \lambda_2}^\ast(h_{\lambda_1} \otimes h_{\lambda_2}) = h_\lambda.$$

Since we have natural $K$-invariant Hilbert structures on every space $H_{\lambda}$, we may dualize $\mu_{\lambda_1, \lambda_2}^\ast$ to obtain nonzero $K$-intertwining maps:
$$\mu_{\lambda_1, \lambda_2} \colon H_\lambda \rightarrow H_{\lambda_1} \otimes H_{\lambda_2}.$$
    Let $\langle \cdot, \cdot \rangle_{\lambda_1 \otimes \lambda_2}$ denote the inner product on $H_{\lambda_1} \otimes H_{\lambda_2}$ given on simple tensors by:
$$\langle h_1 \otimes h_2, h_1' \otimes h_2' \rangle_{\lambda_1 \otimes \lambda_2} = \langle h_1, h_1' \rangle_{\lambda_1} \cdot \langle h_2, h_2' \rangle_{\lambda_2}.$$
Then, for all $h \in H_\lambda$, $\mu_{\lambda_1, \lambda_2}(h)$ satisfies, for all $h_1 \in H_{\lambda_1}$, and $h_2 \in H_{\lambda_2}$:
$$\langle \mu_{\lambda_1, \lambda_2}(h), h_1 \otimes h_2 \rangle_{\lambda_1 \otimes \lambda_2} = \langle h, h_1 \cdot h_2 \rangle_{\lambda_1 + \lambda_2}.$$
In particular, $\mu_{\lambda_1, \lambda_2}(h_\lambda^\vee)$ is a vector in the $\lambda$ weight space of $H_{\lambda_1} \otimes H_{\lambda_2}$ (which has dimension $1$).  Moreover, we compute:
$$\langle \mu_{\lambda_1, \lambda_2}(h_\lambda^\vee), h_{\lambda_1} \otimes h_{\lambda_2} \rangle = \langle h_\lambda^\vee, h_\lambda \rangle = 1.$$
It follows that:
$$\mu_{\lambda_1, \lambda_2}(h_\lambda^\vee) = h_{\lambda_1}^\vee \otimes h_{\lambda_2}^\vee.$$

\begin{thm}
The family of $K$-intertwining maps $\{ \mu_{\lambda_1, \lambda_2} \}$ satisfies the axioms of co-identity, co-commutativity, and co-associativity.
\end{thm}
\begin{proof}
\textbf {Co-Identity: } Note that $H_0 \isom \complex$ consists only of the constant functions, since $X$ is a projective variety and has no global holomorphic functions besides constants.  Hence we have, for all $f \in H_\lambda$, all $C \in H_0$, and $f' \in H_\lambda$:
\begin{eqnarray*}
\langle \mu_{0, \lambda}(f), C \otimes f' \rangle_{0 \otimes \lambda} & = & \langle f, C f' \rangle_{\lambda}, \\
& = & \langle 1 \otimes f, C \otimes f' \rangle_{0 \otimes \lambda}. \\
\end{eqnarray*}
Thus $\mu_{0, \lambda}(f) = 1 \otimes f$.

\textbf{Co-Commutatitivity: }  Let $\lambda = \lambda_1 + \lambda_2$.  For all $f \in H_\lambda$, all $f_1 \in H_{\lambda_1}$, and all $f_2 \in H_{\lambda_2}$, we have:
\begin{eqnarray*}
\langle \mu_{\lambda_1, \lambda_2}(f), f_1 \otimes f_2 \rangle_{\lambda_1 \otimes \lambda_2} & = & \langle f, f_1 f_2 \rangle_\lambda, \\
& = & \langle f, f_2 f_1 \rangle_\lambda, \\
& = & \langle \mu_{\lambda_2, \lambda_1}(f), f_2 \otimes f_1 \rangle_{\lambda_2 \otimes \lambda_1}.
\end{eqnarray*}
It follows that $comm \circ \mu_{\lambda_1, \lambda_2} = \mu_{\lambda_2, \lambda_1}$.

\textbf{Co-Associativity: } Suppose that $f \in H_{\lambda_1 + \lambda_2 + \lambda_3}$, and $f_i \in H_{\lambda_i}$ for $i = 1,2,3$.  Define maps:
\begin{eqnarray*}
\mu_{\lambda_1, (\lambda_2, \lambda_3)} & = & (Id \otimes \mu_{\lambda_2, \lambda_3}) \circ \mu_{\lambda_1, \lambda_2 + \lambda_3}, \\
\mu_{(\lambda_1, \lambda_2), \lambda_3} & = & (\mu_{\lambda_1 + \lambda_2} \otimes Id) \circ \mu_{\lambda_1 + \lambda_2, \lambda_3}.
\end{eqnarray*}

Then we compute:
\begin{eqnarray*}
\langle \mu_{\lambda_1, (\lambda_2, \lambda_3)}(f), f_1 \otimes (f_2 \otimes f_3) \rangle_{\lambda_1 \otimes (\lambda_2 \otimes \lambda_3)} & = & \langle \mu_{\lambda_1, \lambda_2 + \lambda_3}(f), f_1 \otimes f_2 f_3 \rangle_{\lambda_1 \otimes (\lambda_2 + \lambda_3)}, \\
& = & \langle f, f_1 f_2 f_3 \rangle_{\lambda_1 + \lambda_2 + \lambda_3}.
\end{eqnarray*}
A similar computation yields:
$$ \langle \mu_{(\lambda_1, \lambda_2), \lambda_3}(f), (f_1 \otimes f_2) \otimes f_3 \rangle_{(\lambda_1 \otimes \lambda_2) \otimes \lambda_3} = \langle f, f_1 f_2 f_3 \rangle_{\lambda_1 + \lambda_2 + \lambda_3}.$$
It follows that:
$$assoc \circ \mu_{\lambda_1, (\lambda_2, \lambda_3)} = \mu_{(\lambda_1, \lambda_2), \lambda_3}.$$
\end{proof}

The above proposition immediately shows that the ``algebraic modular forms'', considered by Gross in \cite{Gr1}, form a graded ring.  This was also proven, albeit in somewhat different language, but by similar techniques, by Khuri-Makdisi in \cite{KMK}.

\begin{cor}
Suppose that $\variety{G}$ is a connected semisimple algebraic group over $\rats$, and $G = \variety{G}(\reals)$ is compact.  Then for any arithmetic (and hence finite) subgroup $\Gamma \subset G$, the space of modular forms $\M(\Gamma \backslash G)$ of all weights, forms a ring graded by the monoid of dominant weights of $G$.
\end{cor}

\section{Noncompact Groups, Discrete series}

In this section we find a family of comultiplications for discrete series representations of noncompact groups.  However, these intertwining operators arise in the category of unitary representations on Hilbert spaces.  A failure of discrete decomposability prevents these operators from yielding similar operators in the category of $\gK$-modules, except in some very special cases.

Just as the construction for compact groups relied on a geometric realization of representations, we rely on a geometric realization of discrete series representations here.

\subsection{Realization of Discrete Series}
Suppose that $\lambda \in \Lambda_k^+$.  Let $H_\lambda$ denote the irreducible representation of $K$ with highest weight $\lambda$, as defined in the previous section.  Recall that $\tau = \tau_\lambda$ denotes the action of $K$ on $H_\lambda$ by left-translation.  This yields a vector bundle over the symmetric space $Y = G/K$ in the usual way; the smooth sections are defined by:
$$C^\infty(Y, H_\lambda) = \{ f \in C^\infty(G, H_\lambda) \mbox{ such that } \forall k \in K, g \in G, f(gk) = [\tau_\lambda(k^{-1})](f(g))\}.$$

We may define the space of $L^2$ sections similarly, using the previously-defined Hilbert space structure on every $H_\lambda$.  The Killing form yields a measure on $Y = G/K$, normalizing the inner product on $L^2$ sections.  If $\Omega$ denotes the Casimir element of $\ZU$, then define the space:
$$\J(Y, H_\lambda) = \{ f \in [C^\infty \cap L^2](Y, H_\lambda) \mbox{ such that } \Omega f = (\norm{\lambda}^2 - \norm{\rho}^2) f \}.$$

As explained in \cite{O-V}, based on the original work of \cite{H-P} and extended by \cite{Wal}, we have:
\begin{thm}
Suppose that $H_\lambda$ is isomorphic to the lowest K-type of a discrete series representation of $G$.  Then $\J(Y, H_\lambda)$ is a closed subspace of $L^2(Y, H_\lambda)$.  When endowed with the resulting Hilbert space structure and the left-translation action of $G$, $\J(Y, H_\lambda)$ is a irreducible unitary representation of $G$, in the discrete series, with lowest $K$-type isomorphic to $H_\lambda$.
\end{thm}

Hereafter, we define $\J_\lambda$ to be the subspace of $K$-finite vectors in $\J(Y, H_\lambda)$, viewed as an irreducible unitary Harish-Chandra module.

\subsection{Tensor products of discrete series}

Now, suppose that $\lambda_1, \lambda_2 \in \Lambda_k^+$, and $\lambda = \lambda_1 + \lambda_2$.  Theorem 1 of \cite{Var} implies the following:
\begin{prop}
Suppose that $H_{\lambda_1}$, $H_{\lambda_2}$, and $H_\lambda$ are all lowest $K$-types of discrete series representations of $G$.  Then there exists a nonzero $G$-intertwining continuous embedding:
$$\nu_{1,2} \colon \J(Y, H_\lambda) \hookrightarrow \J(Y, H_{\lambda_1}) \hat \otimes \J(Y, H_{\lambda_2}).$$
\end{prop}
\begin{proof}
The completed tensor product $\J(Y, H_{\lambda_1}) \hat \boxtimes \J(Y, H_{\lambda_2})$ is a discrete series representation of $G \times G$, with lowest $(K \times K)$-type $H_{\lambda_1} \boxtimes H_{\lambda_2}$.  Moreover, $H_{\lambda}$ occurs in the restriction of $H_{\lambda_1} \boxtimes H_{\lambda_2}$ from $K \times K$ to $K$, as is well-known, and realized via the map $\mu_{\lambda_1, \lambda_2}$ from before.  Theorem 1 of \cite{Var} applies, and yields the proposition immediately.
\end{proof}

Unfortunately, the results and methods in \cite{Var} do not provide an explicit embedding of these representations.  However, results of \cite{O-V} suggest that such an explicit embedding is possible.

We fix some specialized notation for restriction maps and tensor products.  Suppose that $H_1, H_2, H_3$ are finite-dimensional representations of $K$.  Suppose that $j_1, j_2, j_3$ are sections of the associated bundles:  $j_i \in C^\infty(Y, H_i)$ for $i = 1,2,3$.  Then we may take the fibrewise tensor product; we write $j_1 \odot j_2$ for the element of $C^\infty(Y, H_1 \otimes H_2)$, given by $[j_1 \odot j_2](y) = j_1(y) \otimes j_2(y)$.  We use this notation to distinguish this operation from the ``external'' tensor product $j_1 \otimes j_2 \in C^\infty(Y, H_1) \otimes C^\infty(Y, H_2)$.

There are canonical isomorphsims:
$$comm^Y \colon C^\infty(Y, H_1 \otimes H_2) \rightarrow C^\infty(Y, H_2 \otimes H_1),$$
$$assoc^Y \colon C^\infty(Y, (H_1 \otimes H_2) \otimes H_3) \rightarrow C^\infty(Y, H_1 \otimes (H_2 \otimes H_3)).$$
We have the following properties of commutativity and associativity:
$$comm^Y(j_1 \odot j_2) = j_2 \odot j_1,$$
$$assoc^Y((j_1 \odot j_2) \odot j_3)) = j_1 \odot (j_2 \odot j_3).$$

Then, Theorem 1 of \cite{O-V}, together with the density of the algebraic tensor product in the completed tensor product, directly implies the following:
\begin{prop}
\label{RL2}
The map $j_1 \otimes j_2 \mapsto j_1 \odot j_2$ extends uniquely to a continuous $G$-intertwining linear map of Hilbert spaces:
$$r \colon \J(Y, H_{\lambda_1}) \hat \otimes \J(Y, H_{\lambda_2}) \rightarrow L^2(Y, H_{\lambda_1} \otimes H_{\lambda_2}).$$
\end{prop}
\begin{proof}
The map $j_1 \otimes j_2 \mapsto j_1 \odot j_2$ can be seen as the restriction, to the diagonally embedded $G \subset G \times G$, of an appropriate element of the Hotta-Parthasarathy realization of the discrete series of $G \times G$.  Thus, Theorem 1 of \cite{O-V} applies.
\end{proof}

One has a $G$-intertwining continuous linear map, adjoint to restriction:
$$r^\star \colon L^2(Y, H_{\lambda_1} \otimes H_{\lambda_2}) \rightarrow \J(Y, H_{\lambda_1}) \hat \otimes \J(Y, H_{\lambda_2}).$$

Recall that there are canonical $K$-intertwining maps:
\begin{eqnarray*}
\mu_{\lambda_1,\lambda_2} & \colon & H_\lambda \rightarrow H_{\lambda_1} \otimes H_{\lambda_2}, \\
\mu_{\lambda_1, \lambda_2}^\ast & \colon & H_{\lambda_1} \otimes H_{\lambda_2} \rightarrow H_\lambda.
\end{eqnarray*}
Applying $\mu_{\lambda_1, \lambda_2}$ or $\mu_{\lambda_1, \lambda_2}^\ast$ fibrewise, It follows that there are continuous $G$-intertwining maps:
\begin{eqnarray*}
\mu_{\lambda_1, \lambda_2}^Y & \colon & L^2(Y,H_{\lambda}) \rightarrow L^2(Y, H_{\lambda_1} \otimes H_{\lambda_2}), \\
\mu_{\lambda_1, \lambda_2}^{\ast Y} & \colon & L^2(Y, H_{\lambda_1} \otimes H_{\lambda_2}) \rightarrow  L^2(Y,H_{\lambda}).
\end{eqnarray*}
The maps $\mu_{\lambda_1, \lambda_2}^Y$ are co-commutative and co-associative in the natural way, following the same properties of $\mu_{\lambda_1, \lambda_2}$.

It follows that we may define:
\begin{defn}
There is a continuous $G$-intertwining map:
$$\nu_{\lambda_1, \lambda_2} = r^\star \circ \mu_{\lambda_1, \lambda_2}^Y \colon \J(Y, H_\lambda) \rightarrow \J(Y, H_{\lambda_1}) \hat \otimes \J(Y, H_{\lambda_2}).$$
\end{defn}
If $j \in \J(Y, H_\lambda)$, then $\nu_{\lambda_1, \lambda_2}(j)$ is the unique element of $\J(Y, H_{\lambda_1}) \hat \otimes \J(Y, H_{\lambda_2})$, which satisfies:
$$\langle \nu_{\lambda_1, \lambda_2}(j), j_1 \otimes j_2 \rangle = \langle \mu_{\lambda_1, \lambda_2}^Y j, j_1 \odot j_2 \rangle,$$
for all $j_1 \in \J(Y, H_{\lambda_1})$ and all $j_2 \in \J(Y, H_{\lambda_2})$.

\subsection{Discrete decomposability}

In the category of unitary representations of $G$ on Hilbert spaces, we have constructed a continuous $G$-intertwining map:
$$\nu_{\lambda_1, \lambda_2} \colon \J(Y, H_\lambda) \rightarrow \J(Y, H_{\lambda_1}) \hat \otimes \J(Y, H_{\lambda_2}).$$
However, one may not construct an intertwining map of $\gK$-modules in the na\"ive fashion.  It is true that one may pass the the subspaces of $K$-finite vectors, to get a map:
$$\nu_{\lambda_1, \lambda_2} \colon \J_\lambda \rightarrow \left( \J(Y, H_{\lambda_1}) \hat \otimes \J(Y, H_{\lambda_2}) \right)^K.$$
However, the target of this map rarely coincides with $\J_{\lambda_1} \otimes \J_{\lambda_2}$; this reflects the difference between $K$-finite vectors and $(K \times K)$-finite vectors.  We are very thankful to Gordan Savin for pointing out this difficulty.

The seriousness of this difficulty is illustrated by the following, which easily follows from deeper results of Kobayashi \cite{Kob}:
\begin{prop}
The following conditions are equivalent:
\begin{enumerate}
\item
The unitary representation $\J(Y, H_{\lambda_1}) \hat \otimes \J(Y, H_{\lambda_2})$ of $G$ decomposes discretely, i.e., as a Hilbert space direct sum of irreducible unitary representations of $G$.
\item
The space of $K$-finite vectors coincides with the $K \times K$-finite vectors, in the sense that:
$$\left( \J(Y, H_{\lambda_1}) \hat \otimes \J(Y, H_{\lambda_2}) \right)^K = \J_{\lambda_1} \otimes \J_{\lambda_2}.$$
\item
The intertwining map $\nu_{\lambda_1, \lambda_2}$ restricts to a homomorphism of $\gK$-modules:
$$\J_\lambda \rightarrow \J_{\lambda_1} \otimes \J_{\lambda_2}.$$
\end{enumerate}
\end{prop}
\begin{proof}
To check that (1) $\Rightarrow$ (2), we assume that $\J(Y, H_{\lambda_1}) \hat \otimes \J(Y, H_{\lambda_2})$ decomposes discretely.  Then, using Theorem 4.2 of \cite{Kob}, and the fact that we work with discrete series representations, it follows that the restriction of $\J(Y, H_{\lambda_1}) \hat \otimes \J(Y, H_{\lambda_2})$ to $K$ is admissible.  From this, it follows from Proposition 1.6 of \cite{Kob} that the $K$-finite vectors coincide with the $K \times K$-finite vectors as needed.

The implication (2) $\Rightarrow$ (3) is obvious, following our previous comments.

For the implication (3) $\Rightarrow$ (1), the irreducibility of $\J_\lambda$, the nonvanishing of $\nu_{\lambda_1, \lambda_2}$ (to be proven in Theorem \ref{NV}), and Lemma 1.5 of Kobayashi \cite{Kob} imply that $\left( \J(Y, H_{\lambda_1}) \hat \otimes \J(Y, H_{\lambda_2}) \right)^K$ decomposes discretely as a $\Lie{g}$-module (see Definition 1.1 of \cite{Kob}).  By Lemma 1.3 of \cite{Kob}, this implies a direct sum decomposition:
$$\left( \J(Y, H_{\lambda_1}) \hat \otimes \J(Y, H_{\lambda_2}) \right)^K = \bigoplus_{i=1}^\infty X_i,$$
into orthogonal irreducible $\Lie{g}$-modules.  By Theorem 4.2 of \cite{Kob}, $\left( \J(Y, H_{\lambda_1}) \hat \otimes \J(Y, H_{\lambda_2}) \right)$ is $K$-admissible, and so the Hilbert space $\J(Y, H_{\lambda_1}) \hat \otimes \J(Y, H_{\lambda_2})$ can be identified with the completion of its subspace of $K$-finite vectors.  Thus, we may identify this Hilbert space with the Hilbert space direct sum of completions:
$$\J(Y, H_{\lambda_1}) \hat \otimes \J(Y, H_{\lambda_2}) = \widehat{\bigoplus}_{i=1}^\infty \widehat{X_i}.$$
The subspaces $\widehat{X_i}$ are irreducible unitary representations of $G$, since their underlying $\gK$-modules are irreducible.  Thus (3) $\Rightarrow$ (1).

\end{proof}

The conditions of the previous proposition are quite restrictive.  They are satisfied for holomorphic discrete series representations, for example.  Kobayashi gives many sufficient and practical conditions for discrete decomposability, which we do not recall here.  Instead, we simply give:
\begin{defn}
The pair $\lambda_1, \lambda_2$ is discretely decomposable, if the conditions of the previous proposition are satisfied.
\end{defn}

\subsection{The comultiplication data}

We now prove that the family of comultiplications, given by $\nu_{\lambda_1, \lambda_2}$ is well-behaved:
\begin{prop}
Suppose that $C$ is a submonoid of $\Lambda_k^+$, containing $0$, such that for every $0 \neq \lambda \in C$, $H_\lambda$ is isomorphic to the lowest $K$-type of a discrete series representation of $G$.  Suppose that every pair $\lambda_1, \lambda_2 \in C$ is discretely decomposable.  Then the family of comultiplications $\{ \nu_{\lambda_1, \lambda_2} \}$, among the family of Harish-Chandra modules $\J_\lambda$ ($\lambda \in C$) has the properties of co-identity, co-commutativity, and co-associativity.
\end{prop}
\begin{proof}
\textbf{Co-Identity: } We always fix $\J_0 = \complex$, and we define the comultiplications $\nu_{0, \lambda}$ to be the canonical maps.  Thus the property of co-identity is vacuous.

\textbf{Co-Commutatitivity: }  Let $\lambda = \lambda_1 + \lambda_2$, with $\lambda_1, \lambda_2 \in C$.  Fix $j \in \J_\lambda$, $j_1 \in \J_{\lambda_1}$ and $j_2 \in \J_{\lambda_2}$.  Then we have:
\begin{eqnarray*}
\langle \nu_{\lambda_1, \lambda_2}(j), j_1 \otimes j_2 \rangle,
& = & \langle \mu_{\lambda_1, \lambda_2}^Y j, j_1 \odot j_2 \rangle, \\
& = & \langle comm^Y(\mu_{\lambda_1, \lambda_2}^Y j), comm^Y(j_1 \odot j_2) \rangle, \\
& = & \langle \mu_{\lambda_2, \lambda_1}^Y j, j_2 \odot j_1 \rangle,  \\
& = & \langle \nu_{\lambda_2, \lambda_1}(j), j_2 \otimes j_1 \rangle.
\end{eqnarray*}
Note that the middle step essentially relies on the co-commutativity of $\mu_{\lambda_1, \lambda_2}$.  It follows that the family $\{ \nu_{\lambda_1, \lambda_2} \}$ is co-commutative.

\textbf{Co-Associativity: } Suppose that $\lambda_1, \lambda_2, \lambda_3 \in C$.  Suppose that $j \in \J_{\lambda_1 + \lambda_2 + \lambda_3}$, and $j_i \in \J_{\lambda_i}$ for $i = 1,2,3$.  Define maps:
\begin{eqnarray*}
\nu_{\lambda_1, (\lambda_2, \lambda_3)} & = & (Id \otimes \nu_{\lambda_2, \lambda_3}) \circ \nu_{\lambda_1, \lambda_2 + \lambda_3}, \\
\nu_{(\lambda_1, \lambda_2), \lambda_3} & = & (\nu_{\lambda_1,\lambda_2} \otimes Id) \circ \nu_{\lambda_1 + \lambda_2, \lambda_3}.
\end{eqnarray*}

Then we compute:
\begin{eqnarray*}
\langle \nu_{\lambda_1, (\lambda_2, \lambda_3)}(j), j_1 \otimes (j_2 \otimes j_3) \rangle
& = & \left\langle \left[ Id \otimes \mu_{\lambda_2, \lambda_3}^Y \right] \left( \nu_{\lambda_1, \lambda_2 + \lambda_3} j \right) , j_1 \otimes (j_2 \odot j_3) \right\rangle, \\
& = & \left\langle \nu_{\lambda_1, \lambda_2 + \lambda_3} j, j_1 \otimes \mu_{\lambda_2, \lambda_3}^{\ast Y} (j_2 \odot j_3) \right\rangle, \\
& = & \left\langle \mu_{\lambda_1, \lambda_2+\lambda_3}^Y j, j_1 \odot \mu_{\lambda_2, \lambda_3}^{\ast Y} (j_2 \odot j_3) \right\rangle, \\
& = & \left\langle j, \mu_{\lambda_1, \lambda_2+\lambda_3}^{\ast Y} \left( j_1 \odot \mu_{\lambda_2, \lambda_3}^{\ast Y} (j_2 \odot j_3) \right) \right\rangle, \\
& = & \left\langle j, \mu_{\lambda_1 + \lambda_2, \lambda_3}^{\ast Y} \left( \mu_{\lambda_1, \lambda_2}^{\ast Y} (j_1 \odot j_2) \odot j_3 \right) \right\rangle, \\
& = & \langle \nu_{(\lambda_1, \lambda_2), \lambda_3}(j), (j_1 \otimes j_2) \otimes j_3 \rangle.
\end{eqnarray*}
 We are using the associativity of the bilinear operations $\mu^{\ast Y}$, which follows from ordinary associativity of multiplication, and the Borel-Weil realization.  It follows that the family $\{ \nu_{\lambda_1, \lambda_2} \}$ is co-associative.
\end{proof}

The above proposition yields the following:
\begin{prop}
Let $C$ be any submonoid of $\Lambda_k^+$ containing $0$, such that for all $0 \neq \lambda \in C$, the representation $H_\lambda$ is isomorphic to the lowest $K$-type of some discrete series representation of $G$.  Suppose that every pair in $C$ is discretely decomposable.  Let $\M_C = \bigoplus_{\lambda \in C} \M(\J_\lambda)$.  Then the family of comultiplications $\{ \nu_{\lambda_1, \lambda_2} \}$ endow $\M_C$ with the structure of a commutative ring.
\end{prop}

Of course, for the ring $\M_C$ to be interesting, we should expect the multiplication maps to be non-trivial.  In other words, we must prove that the comultiplications $\nu_{\lambda_1, \lambda_2}$ are nonzero.

\subsection{Nonvanishing of intertwining operators}

In order to prove that the maps $\nu_{\lambda_1, \lambda_2}$ are nonzero, we use spherical trace functions.  Essentially all of the facts we use about spherical trace functions can be found in the excellent treatment of \cite{W-W}.

First, we fix some notation.  Let $\ii \colon H_{\lambda} \hookrightarrow \J(Y, H_\lambda)$ denote the unique $K$-equivariant inclusion, which is adjoint to the evaluation map (at the identity) $\ev \colon \J(Y, H_\lambda) \rightarrow H_{\lambda}$.  Let $\pp \colon \J(Y, H_\lambda) \rightarrow \ii(H_\lambda)$ denote the orthogonal projection onto the lowest $K$-type.

Similarly, let
$$\II \colon H_{\lambda_1} \boxtimes H_{\lambda_2} \rightarrow \J(Y \times Y, H_{\lambda_1} \boxtimes H_{\lambda_2})$$
denote the $K \times K$-equivariant inclusion adjoint to the evaluation map.  Let
$$\PP \colon \J(Y \times Y, H_{\lambda_1} \boxtimes H_{\lambda_2}) \rightarrow \II \left( H_{\lambda_1} \boxtimes H_{\lambda_2} \right)$$
denote the projection onto the lowest $(K \times K)$-type.

There is a close relationship between projecting onto the lowest $K$-type (or $(K \times K)$-type), and evaluating at the identity.  Namely, for all $j \in \J(Y, H_\lambda)$, we have:
$$[\ii^{-1} \circ \pp](j) = {1 \over \dd} \cdot j(1),$$
where $\dd$ is the formal degree of $\J(Y, H_\lambda)$.  A corresponding formula holds for $J \in \J(Y \times Y, H_{\lambda_1} \boxtimes H_{\lambda_2})$:
$$[\II^{-1} \circ \PP](J) = {1 \over \DD} \cdot J(1).$$

Let $\Ell$ denote the left-translation action of $G$ on $\J(Y, H_\lambda)$.  Let $\LL$ denote the left-translation action of $G \times G$ on $\J(Y \times Y, H_{\lambda_1} \boxtimes H_{\lambda_2})$.

Let $h^1, \ldots, h^a$ denote an orthonormal basis of $H_\lambda$.  Let $e^s = \mu_{\lambda_1, \lambda_2}(h^s)$ for all $1 \leq s \leq a$.  Let $F$ denote the orthogonal complement of $\mu_{\lambda_1, \lambda_2}(H_\lambda)$ in $H_{\lambda_1} \otimes H_{\lambda_2}$.  Let $f^1, \ldots, f^b$ denote an orthonormal basis of $F$.  Thus the set $\{e^1, \ldots, e^a, f^1, \ldots, f^b \}$ is an orthogonal basis of $H_{\lambda_1} \otimes H_{\lambda_2}$.  For all $1 \leq s \leq a$, we have $\langle e^s, e^s \rangle = c$, for some fixed positive number $c$.

With these choices, we may write the spherical trace functions quite explicitly.  On $G$, the spherical trace function for $\J(Y, H_\lambda)$ is given by:
\begin{eqnarray*}
\phi(g) & = & Tr(\pp \circ \Ell(g) \circ \pp), \\
& = & Tr(\ii^{-1} \circ \pp \circ \Ell(g) \circ \pp \circ \ii), \\
& = & \sum_{s=1}^a \langle \ii^{-1} \circ \pp \circ \Ell(g) \circ \pp \circ \ii h^s, h^s \rangle, \\
& = & {1 \over \dd} \cdot \sum_{s=1}^a \langle [\Ell(g) \circ \pp \circ \ii h^s](1), h^s \rangle, \\
& = & {1 \over \dd} \cdot \sum_{s=1}^a \langle [\pp \circ \ii h^s](g^{-1}), h^s \rangle, \\
& = & {1 \over \dd} \cdot \sum_{s=1}^a \langle [\ii h^s](g^{-1}), h^s \rangle.
\end{eqnarray*}
Similarly, on $G \times G$, we have an explicit expresssion for the spherical trace function.  Restricting to the diagonal $G \subset G \times G$ yields:
\begin{eqnarray*}
\Phi(g) & = & Tr \left( \PP \circ \LL(g) \circ \PP \right), \\
& = & {1 \over {c \DD}} \cdot \left( \sum_{s=1}^a \langle \II e^s(g^{-1}), e^s \rangle + \sum_{t=1}^b \langle \II f^t(g^{-1}), f^t \rangle \right).
\end{eqnarray*}

We recall Lemma 2 of \cite{Var}
\begin{lem}
\label{VL}
The product $\phi(g) \cdot \Phi(g)$ is an integrable function on $G$ (with respect to a Haar measure).  Moreover,
$$0 < \int_G \phi(g) \overline{\Phi(g)} dg \in \reals.$$
\end{lem}
\begin{proof}
This is precisely the statement of Lemma 2 of \cite{Var}, for the pair of groups $G \subset G \times G$.
\end{proof}

Using Schur's orthogonality relations, we may prove:
\begin{lem}
For all $1 \leq s \leq a$, and for all $1 \leq t \leq b$, we have:
$$\int_G \langle [\ii h^s](g^{-1}), h^s \rangle \overline{ \langle \II f^t(g^{-1}), f^t \rangle} dg = 0.$$
\end{lem}
\begin{proof}
Using the decomposition $G = KA^+ K$, we can express the above integral as:
$$\int_K \int_{A^+} \Delta(a) \int_K \langle [\ii h^s](k_1 a k_2)^{-1}, h^s \rangle \overline{ \langle \II f^t(k_1 a k_2)^{-1}, f^t \rangle} dk_2 da dk_1.$$
Consider the innermost integral above:
$$\int_K \langle [\Ell(k_2)[\ii h^s]](k_1 a)^{-1}, h^s \rangle \overline{ \langle [\LL(k_2) [\II f^t]](k_1 a)^{-1}, f^t \rangle} dk_2.$$
Since $\ii h^s$, and $\II f^t$ are lowest $K$-type and $(K \times K)$-type vectors, respectively, the integral above equals:
$$\int_K \langle \tau_\lambda(k_2) h, h^s \rangle \overline{ \langle (\tau_{\lambda_1} \otimes \tau_{\lambda_2})(k_2) f, f^t \rangle } dk_2,$$
where $h = [\ii h^s](k_1 a)^{-1}$ and $f$ is the orthogonal projection of $[\II f^t](k_1 a)^{-1}$ in $F$; note that $h \in H_\lambda$ and $f \in F$.

The first inner product is a matrix coefficient for the representation $H_\lambda$ of $K$.  The second inner product is a matrix coefficient for the representation $F$ of $K$.  Since $H_\lambda$ occurs with multiplicity one in $H_{\lambda_1} \otimes H_{\lambda_2}$, and $F$ is complementary to the embedding $\mu_{\lambda_1, \lambda_2}(H_\lambda)$, these matrix coefficients are orthogonal by Schur's orthogonality.  The lemma follows immediately.
\end{proof}

Again, using Schur's orthogonality relations, we prove:
\begin{lem}
For all $1 \leq s,s' \leq a$, we have:
$$\int_G \langle [\ii h^s](g^{-1}), h^s \rangle \overline{ \langle [\II e^{s'}](g^{-1}), e^{s'} \rangle} dg = 0, \mbox{ if } s \neq s'.$$
If $1 \leq s \leq a$, we have:
$$\int_G \langle [\ii h^s](g^{-1}), h^s \rangle \overline{ \langle [\II e^s](g^{-1}), e^s \rangle} dg = \int_G \langle \mu_{\lambda_1, \lambda_2} \left( [\ii h^s](g^{-1}) \right), [\II e^s](g^{-1}) \rangle.$$
\end{lem}
\begin{proof}
We use the decomposition $G = KA^+ K$ as before, and compute the above integral as:
$$\int_K \int_{A^+} \Delta(a) \int_K \langle [\ii h^s](k_1 a k_2)^{-1}, h^s \rangle \overline{ \langle [\II e^{s'}](k_1 a k_2)^{-1}, e^{s'} \rangle} dk_2 da dk_1.$$
The innermost integral is then:
$$\int_K \langle [\Ell(k_2)[\ii h^s]](k_1 a)^{-1}, h^s \rangle \overline{ \langle [\LL(k_2) [\II e^{s'}]](k_1 a)^{-1}, e^{s'} \rangle } dk_2.$$
By Schur's orthogonality relations, using the embedding $\mu_{\lambda_1,\lambda_2}$, the above integral equals:
$$\left\langle \mu_{\lambda_1, \lambda_2} \left( [\ii h^s](k_1 a)^{-1} \right) , [\II e^{s'}](k_1 a)^{-1} \right\rangle \overline{\langle  e^s , e^{s'} \rangle}.$$
Since the inner product is $K$-invariant, the above equals:
$$ \overline{\langle  e^s , e^{s'} \rangle} \cdot \int_K \langle \mu_{\lambda_1, \lambda_2} \left( [\ii h^s](k_1 a k_2)^{-1} \right) , [\II e^{s'}](k_1 a k_2)^{-1} \rangle dk_2.$$

From this, we compute the whole integral:
\begin{eqnarray*}
\int_G \langle [\ii h^s](g^{-1}), h^s \rangle \overline{ \langle [\II e^{s'}](g^{-1}), e^{s'} \rangle} dg \\ = \overline{\langle e^s, e^{s'} \rangle } \cdot \int_G \langle \mu_{\lambda_1, \lambda_2} \left( [\ii h^s](g^{-1}) \right), [\II e^{s'}](g^{-1}) \rangle dg.
\end{eqnarray*}

The lemma follows immediately.
\end{proof}

The previous three lemmas allow us to prove:
\begin{thm}
The map $\nu_{\lambda_1, \lambda_2}$ is nonzero.
\label{NV}
\end{thm}
\begin{proof}
By Lemma \ref{VL}, we have:
$$\int_G \phi(g) \overline{\Phi(g)} dg \neq 0.$$
Expanding this, and disregarding the formal degree constants, yields:
\begin{eqnarray*}
0 & \neq & \left( \sum_{s,s'} \int_G \langle [\ii h^s](g^{-1}), h^s \rangle \overline{ \langle [\II e^{s'}](g^{-1}), e^{s'} \rangle } \right) dg \\
& & + \left( \sum_{s, t} \int_G \langle [\ii h^s](g^{-1}), h^s \rangle \overline{ \langle [\II f^t](g^{-1}), f^t \rangle } \right) dg.
\end{eqnarray*}

The previous two lemmas imply the vanishing of all but a few terms above.  We are left with:
$$0 \neq \sum_s \int_G \langle [\ii h^s](g^{-1}), h^s \rangle \overline{ \langle [\II e^{s}](g^{-1}), e^{s} \rangle }dg .$$
Finally, the previous lemma implies that there exists an $s$ such that:
$$0 \neq \int_G \langle \mu_{\lambda_1, \lambda_2} \left( [\ii h^s](g^{-1}) \right), [\II e^s](g^{-1}) \rangle dg .$$
Hence, if $j = \ii h^s \in \J(Y, H_\lambda)$, and $J = \II e^s \in \J(Y \times Y, H_{\lambda_1} \boxtimes H_{\lambda_2})$, then:
$$\langle \nu_{\lambda_1, \lambda_2} j, J \rangle \neq 0.$$
In particular, $\nu_{\lambda_1, \lambda_2}(j) \neq 0$.
\end{proof}

In fact, the final steps in this proof illustrate something a bit stronger:
\begin{cor}
Suppose that $j$ is a nonzero lowest $K$-type vector in $\J(Y, H_\lambda)$.  Then the projection of $\nu_{\lambda_1, \lambda_2}(j)$ onto the lowest $(K \times K)$-type of $\J(Y \times Y, H_{\lambda_1} \boxtimes H_{\lambda_2})$ is nonzero.
\end{cor}
\begin{proof}
The proof of the previous theorem illustrates that the result is true for some element $\ii h^s$ of the lowest $K$-type of $\J(Y, H_\lambda)$.  Since $\nu_{\lambda_1, \lambda_2}$ is $K$-intertwining, and $H_\lambda$ is irreducible, the result is true for every nonzero lowest $K$-type vector.
\end{proof}

\section{Classical Modular Forms}

With $G$, $\Gamma$, and a discrete series Harish-Chandra module $\J_\lambda$ as before, we have defined a space of modular forms $\M(\J_\lambda)$.  In this section, we discuss these modular forms in a more classical language.

\subsection{Transition to classical modular forms}
\begin{defn}
The space of {\em classical modular forms} of weight $H_\lambda$, denoted $\M_{cl}(H_\lambda)$ is the subspace of functions $f \in  C^\infty(\Gamma \backslash G, H_\lambda)$, satisfying the additional properties:
\begin{itemize}
\item
$f(gk) = \tau(k^{-1})(f(g))$ for all $g \in G$, $k \in K$.
\item
$\Omega f = ( \norm{\lambda}^2 - \norm{\rho}^2 ) f.$
\item
$f$ has moderate growth.
\end{itemize}
\end{defn}

Then we can view our modular forms as classical modular forms:
\begin{prop}
For all $F \in \M(\J_\lambda)$, there exists a unique classical modular form $f \in \M_{cl}(H_\lambda)$, which satisfies:
$$\langle f(g), h \rangle = F(\ii(h))(g), \mbox{ for all } h \in H_\lambda, g \in G.$$
This yields an injective linear map $\M(\J_\lambda) \hookrightarrow \M_{cl}(H_\lambda)$.
\end{prop}
\begin{proof}
Since $H_\lambda$ is a Hilbert space, the values of $f(g)$ are completely determined by the formula $\langle f(g), h \rangle = F(\ii(h))(g)$.  Fixing $F \in \M(\J_\lambda)$, we may compute:
\begin{eqnarray*}
\langle f(gk), h \rangle & = & F(\ii(h))(gk), \\
& = & F(\Ell(k) \ii(h))(g), \\
& = & F(\ii(\tau(k) h))(g), \\
& = & \langle f(g), \tau(k) h \rangle, \\
& = & \langle \tau(k^{-1}) f(g), h \rangle.
\end{eqnarray*}

Next, if $Z \in \ZU$, and $\omega$ is the infinitesimal character of $\J_\lambda$, then we compute:
$$\langle Z f(g), h \rangle = Z[F(\ii(h))](g) = \omega(Z) \cdot F(i(h))(g).$$
In particular, this holds for $\Omega \in \ZU$, so that $f$ is an eigenfunction of $\Omega$ with the same eigenvalue as $F(\ii(h))$.

Finally, $f$ has moderate growth, since fixing an orthonormal basis $h_1, \ldots, h_a$ of $H_\lambda$, we compute:
\begin{eqnarray*}
\vert f(g) \vert_\lambda^2 & = & \sum_{s=1}^a \vert \langle f(g), h_s \rangle \vert^2, \\
& = & \sum_{s=1}^a \vert F(\ii(h_s))(g) \vert^2, \\
& \leq & a C \Vert g \Vert^{2n},
\end{eqnarray*}
where $C$ and $n$ are the suprema of the constants of moderate growth for $F(\ii(h_s))$.

Thus we have a linear map from $\M(\J_\lambda)$ to $\M_{cl}(H_\lambda)$ as desired.  It is injective, since the classical modular form $f$ determines the restriction of a modular form $F$ to its lowest $K$-type, and $\J_\lambda$ is generated as a $\U$-module by its lowest $K$-type.
\end{proof}

\subsection{Holomorphic Modular Forms}

Suppose that $G/K$ is a Hermitian symmetric domain.  In this case, multiplication of scalar-valued holomorphic modular forms is well-understood.  Representation theoretically, this multiplication can be associated to properties of the holomorphic discrete series of $G$.

Let $\Lie{g} = \Lie{k} \oplus \Lie{p}$ denote the Cartan decomposition.  There exists a decomposition of the roots of $(\Lie{g}, \Lie{t})$ into positive and negative roots, in such a way that the sum of any two noncompact positive roots is not a root.  This yields a natural decomposition $\Lie{p} = \Lie{p}^+ \oplus \Lie{p}^-$.  The elements of $\Lie{p}^-$ are called the ``lowering operators''.  Holomorphic functions on $G/K$ are precisely those functions which are annihilated by all lowering operators.

The scalar-valued holomorphic discrete series is a family of discrete series representations $\J_w$ of $G$, parameterized by the nonzero elements $w \in C \subset \nats$ in a submonoid of $\nats$, which satisfy the following properties:
\begin{itemize}
\item
The lowest $K$-type of $\J_w$ is one-dimensional, denoted $\complex_w$.  The action of $K$ on $\complex_w$ will be denoted $\chi_w$.
\item
Every element $j \in \complex_w$, the lowest $K$-type, is annihilated by the lowering operators.
\item
Every pair $w_1, w_2 \in C$ is discretely decomposable (cf. \cite{Rep}).
\end{itemize}

For consistency, we realize the holomorphic discrete series, as all discrete series representations, as sections of a bundle over $G/K$, which satisfy a growth condition and the second-order differential equation from the Casimir operator.  In this realization, we have the ``evaluation at the identity'' operator $\ev_w \colon \J_w \rightarrow \complex_w$, and the adjoint inclusion $\ii_w \colon \complex_w \rightarrow \J_w$.  Define $j_w^\vee = \ii_w(1)$; then
$$j_w^\vee(1) = \langle j_w^\vee, j_w^\vee \rangle = \delta_w,$$
where $\delta_w$ is the formal degree of $\J_w$.  Let $j_w = {1 \over {\delta_w}} j_w^\vee$, so that $j_w$ is the unique lowest $K$-type vector satisfying
$$j_w(1) = \langle j_w, j_w^\vee \rangle = 1.$$

When $w_1 + w_2 = w$, and $w_1, w_2 \in C$, we have constructed intertwining operators
$$\nu_{w_1, w_2} \colon \J_w \rightarrow \J_{w_1} \otimes \J_{w_2}.$$
We prove the following fact about these intertwining operators:
\begin{lem}
Let $\nu_{w_1,w_2}$ be the canonical intertwining operator from $\J_w$ to $\J_{w_1} \otimes \J_{w_2}$.  Then
$$\nu_{1,2}(j_w^\vee) = j_{w_1}^\vee \otimes j_{w_2}^\vee.$$
\end{lem}
\begin{proof}
We first note that $j_{w_1} \odot j_{w_2}$ is a holomorphic element of $C^\infty(Y, \complex_w)$, since holomorphy is preserved under products.  As such, it is annihilated by the Dirac operator, and hence is an eigenfunction of the Casimir operator, of the appropriate eigenvalue for $\J_w$ by Proposition 3.2 of Parthasarathy \cite{Par}.  It is also in $L^2(Y, \complex_w)$, by Proposition \ref{RL2}.  Finally, $(j_{w_1} \odot j_{w_2})(1) = 1$, by simply multiplying.  Hence $j_{w_1} \odot j_{w_2}$ is a normalized lowest $K$-type vector in $\J_w$, and from uniqueness of such, it follows that:
$$j_{w_1} \odot j_{w_2} = j_w.$$

To prove the adjoint statement, we note that $\nu_{1,2}(j_w^\vee)$ is a lowest $K$-type vector in $\J_{w_1} \otimes \J_{w_2}$, since the appropriate $K$-type has multiplicity one (as noted by Repka \cite{Rep}, for example).  It follows that for some nonzero constant $C$, we have:
$$\nu_{1,2}(j_w^\vee) = C j_1^\vee \otimes j_2^\vee.$$
To prove that this constant equals one, we compute:
\begin{eqnarray*}
\langle \nu_{1,2}(j_w^\vee), j_{w_1} \otimes j_{w_2} \rangle & = & \langle \mu^Y(j_w^\vee), j_{w_1} \odot j_{w_2} \rangle, \\
& = & \langle j_w^\vee, j_w \rangle, \\
& = & 1.
\end{eqnarray*}
On the other hand, we compute:
\begin{eqnarray*}
\langle \nu_{1,2}(j_w^\vee), j_{w_1} \otimes j_{w_2} \rangle & = & \langle C j_{w_1}^\vee \otimes j_{w_2}^\vee, j_{w_1} \otimes j_{w_2} \rangle, \\
& = & C \langle j_{w_1}^\vee, j_{w_1} \rangle \langle j_{w_2}^\vee, j_{w_2} \rangle, \\
& = & C.
\end{eqnarray*}
It follows that $C = 1$.
\end{proof}

When working with holomorphic discrete series, we may consider modular forms that satisfy a first-order differential equation, instead of a second-order equation.
\begin{defn}
The space of holomorphic modular forms of weight $w$, denoted $\M_{hol}(w)$ is the space of smooth functions $f \colon \Gamma \backslash G \rightarrow \complex$ such that:
\begin{itemize}
\item
$f(gk) = \chi_w(k^{-1})(f(g))$ for all $g \in G$, $k \in K$.
\item
$X f = 0$, for all $X \in \Lie{p}^-$.
\item
$f$ has moderate growth.
\end{itemize}
\end{defn}

Holomorphic modular forms form a ring.  Indeed, it is a fortunate fact that if $f_1 \in \M_{hol}(w_1)$ and $f_2 \in \M_{hol}(w_2)$, then $f_1 \cdot f_2 \in \M_{hol}(w_1 + w_2)$; not only is the pointwise product naturally a section of the bundle $\complex_w$, but it also satisfies the appropriate differential equations for holomorphy.  Thus, multiplication of holomorphic modular forms is ``easy''.

Again, by computations in \cite{Par}, elements $f \in \M_{hol}(w)$ also satisfy $\Omega f = ( \norm{\lambda}^2 - \norm{\rho}^2 ) f,$ for some $\lambda$ corresponding to $w$.  It follows that $\M_{hol}(w) \subset \M_{cl}(\complex_w)$, for all $w \in C$.  Moreover, we have:
\begin{prop}
Suppose that $F \in \M(\J_w)$ is a modular form, where $\J_w$ is the scalar-valued holomorphic discrete series with lowest $K$-type $\complex_w$.  Let $f \in \M_{cl}(\complex_w)$ be the associated classical modular form.  Then $f \in \M_{hol}(w)$.
\end{prop}
\begin{proof}
Every lowest $K$-type vector in $\J_w$ is annihilated by the ``lowering operators'' $X \in \Lie{p}^-$.  If $F \in \M(\J_w)$ is a modular form, then the associated classical modular form is given by a simple formula:
$$f(g) = F(\ii(1))(g) = F(j_w^\vee)(g).$$
It follows that:
$$X f(g) = X F(j_w^\vee)(g) = F( X j_w^\vee)(g) = 0,$$
for every $X \in \Lie{p}^-$.
\end{proof}

It is not clear that ``easy multiplication'' of holomorphic modular forms is compatible with the multiplication induced from the intertwining operator $\nu_{w_1, w_2}$.  For $w_1, w_2 \in C$, and $w = w_1 + w_2$, we have two bilinear maps:
\begin{eqnarray*}
m_{hol} & \colon & \M_{hol}(w_1) \otimes \M_{hol}(w_2) \rightarrow \M_{hol}(w), \\
m_{1,2} & \colon & \M(\J_{w_1}) \otimes \M(\J_{w_2}) \rightarrow \M(\J_{w}).
\end{eqnarray*}
We also have linear maps, for all $w \in C$, given by ``classicalization'':
$$cl_w \colon \M(\J_w) \hookrightarrow \M_{hol}(w).$$
We now prove compatibility:
\begin{prop}
The following diagram commutes, for any pair $w_1, w_2 \in C$:
$$
\xymatrix@C+1.5cm{
  \M(\J_{w_1}) \otimes \M(\J_{w_2}) \ar[d]_{cl_{w_1} \otimes cl_{w_2}} \ar[r]^{m_{1,2}}
                & \M(\J_w) \ar[d]^{cl_w}  \\
  \M_{hol}(w_1) \otimes \M_{hol}(w_2) \ar[r]^{m_{hol}}
                & \M_{hol}(w)             }
$$
\end{prop}
\begin{proof}
Suppose that $F_1 \in \M(\J_{w_1})$ and $F_2 \in \M(\J_{w_2})$.  Let $f_1, f_2$ be the corresponding holomorphic classical modular forms.  Then we have:
$$f_1(g) = F_1(j_{w_1}^\vee)(g), \mbox{ and } f_2(g) = F_2(j_{w_2}^\vee)(g).$$
Hence, we have:
$$[f_1 \cdot f_2](g) = F_1(j_{w_1}^\vee)(g) \cdot F_2(j_{w_2}^\vee)(g).$$
On the other hand, let $F$ be the product of $F_1$ and $F_2$ via $m_{1,2}$.  Let $f$ be the corresponding holomorphic classical modular form.  Then we have:
\begin{eqnarray*}
f(g) & = & F(j_w^\vee)(g), \\
& = & m( [F_1 \otimes F_2] \nu_{w_1,w_2}(j_w^\vee) )(g), \\
& = & m( [F_1 \otimes F_2] (j_{w_1}^\vee \otimes j_{w_2}^\vee ) )(g), \\
& = & F_1(j_{w_1}^\vee)(g) \cdot F_2(j_{w_2}^\vee)(g).
\end{eqnarray*}
\end{proof}

Therefore, our multiplication of modular forms, via intertwining operators $\nu_{w_1, w_2}$, agrees with ordinary multiplication of holomorphic modular forms.

When $G/K$ is a rational tube domain, the quotient $Y(\Gamma) = \Gamma \backslash G / K$ is a locally symmetric space.  Its Satake compactification may be defined, as a projective algebraic variety, by:
$$\overline{Y(\Gamma)} = Proj(\M) = Proj \bigoplus_{w = 0}^\infty \M(\J_w),$$
where $\M$ is the ring of scalar-valued holomorphic modular forms.  Thus the Satake compactification may be constructed completely within the framework of representation theory.  We are very interested in whether the geometric properties of $\overline{Y(\Gamma)}$ may be recovered via representation theory as well.  For example, it can be shown that the formal degree $\delta_w$ of $\J_w$ satisfies $\delta_w = O(w^n)$, where $n = dim(G/K)$.  Heuristically, $dim(\M(\J_w)) = O(\delta_w)$ -- this is precisely true when considering cusp forms, or $L^2$ forms, in the limit as $\Gamma$ becomes smaller.  Therefore, the graded ring $\M$ has dimension $n$, and so we may accurately predict that:
$$dim(\overline{Y(\Gamma)}) = dim(G/K).$$

\subsection{Multiplication of Vector-Valued Modular Forms}

Consider, as previously, the case when $G/K$ is a rational tube domain.  In addition to scalar-valued modular forms, we may consider vector-valued modular forms.
\begin{defn}
Suppose that $H_\lambda$ is an irreducible unitary representation of $K$.  The space of holomorphic modular forms of weight $H_\lambda$, denoted $\M_{hol}(H_\lambda)$ is the space of smooth functions $f \colon \Gamma \backslash G \rightarrow H_\lambda$ such that:
\begin{itemize}
\item
$f(gk) = \tau_\lambda(k^{-1})(f(g))$ for all $g \in G$, $k \in K$.
\item
$X f = 0$, for all $X \in \Lie{p}^-$.
\item
$f$ has moderate growth.
\end{itemize}
\end{defn}

As before, if $H_\lambda$ is the lowest $K$-type of a holomorphic discrete series representation $\J_\lambda$, then there is an injective classicalization map:
$$\M(\J_\lambda) \rightarrow \M_{hol}(H_\lambda).$$

As tensor products of holomorphic discrete series representations are discretely decomposable, one arrives at a ring of modular forms using intertwining operators among the $\J_\lambda$.  Classically, multiplication of vector-valued modular forms is not as well-understood as for scalar-valued modular forms.  However, our ring structure suggests (is compatible with) the following multiplication:
\begin{defn}
Suppose that $f_1 \in \M_{hol}(H_{\lambda_1})$ and $f_2 \in \M_{hol}(H_{\lambda_2})$ are two holomorphic modular forms; thus $H_{\lambda_1}$, $H_{\lambda_2}$, and $H_\lambda$ ($\lambda = \lambda_1 + \lambda_2$) are lowest $K$-types of discrete series representations.  Define the product by:
$$[f_1 \cdot f_2](g) = \mu_{\lambda_1, \lambda_2}^\ast (f_1(g) \otimes f_2(g)) \in H_{\lambda}.$$
\end{defn}

Here, we recall that $\mu_{\lambda_1, \lambda_2}$ is the canonical intertwining operator obtained via the Borel-Weil realization.  The resulting function is in $\M_{hol}(H_\lambda)$, since holomorphicity is preserved.

The classical ring structure, defined above to correspond with our general theory, seems to be new.  One arrives at a ring
$$\M^{vec} = \bigoplus_{\lambda \in \Lambda_k^+} \M_{hol}(H_\lambda),$$
graded by the abelian group $\Lambda_k = \Lambda$.  Viewing elements $\lambda \in \Lambda$ as algebraic characters $T_\complex \rightarrow \complex^\times$, we have an action of $T_\complex$ on the ring $\M$.  Namely, if $f \in \M_{hol}(H_\lambda)$, and $t \in T_\complex$, we define:
$$t[f] = \lambda(t) f.$$

This action yields a geometric action.  The complex algebraic torus $T_\complex$ acts on the $\complex$-scheme $Y^{vec} = Spec(\M^{vec})$.  It seems important to consider the following questions:
\begin{itemize}
\item
Describe the geometry of the GIT quotient $Y^{vec} // T_\complex$.  Does it provide a compactification of the locally symmetric space $\Gamma \backslash G / K$?
\item
Does this quotient have a natural model over a finite extension of $\rats$?  Does it coincide with a toroidal compactification?
\item
Describe the behavior of constant terms of vector-valued modular forms under multiplication.  How does this reflect the cuspidal geometry of $Y^{vec} // T_\complex$?
\end{itemize}

We hope to approach these questions further in a future paper.

\bibliographystyle{amsalpha}
\bibliography{MultMod}
\end{document}